\documentstyle[12pt,twoside]{article}
\font\teneufm=eufm10 scaled \magstep1
\font\seveneufm=eufm7 scaled \magstep1
\font\fiveeufm=eufm5  scaled \magstep1
\newfam\eufmfam
\textfont\eufmfam=\teneufm
\scriptfont\eufmfam=\seveneufm
\scriptscriptfont\eufmfam=\fiveeufm
\def\frak#1{{\fam\eufmfam\relax#1}}
\newfam\msbfam
\font\tenmsb=msbm10 scaled \magstep1  \textfont\msbfam=\tenmsb
\font\sevenmsb=msbm7 scaled \magstep1 \scriptfont\msbfam=\sevenmsb
\font\fivemsb=msbm5 scaled \magstep1  \scriptscriptfont\msbfam=\fivemsb
\def\Bbb{\fam\msbfam \tenmsb}

\def\RR{{\Bbb R}}
\def\CC{{\Bbb C}}

\def\HollowBoxx #1#2#3{{\dimen0=#1 \advance\dimen0 by -#2
       \dimen1=#1 \advance\dimen1 by #3
        \vrule height 0pt depth #3 width #2
       \hskip -#3
       \vrule height #1 depth #3 width #3}}
 \def\LeftContraction{\mathord{\kern1.45pt \HollowBoxx{6pt}{3.5pt}{.4pt}}\,}
 \def\HollowBox #1#2#3{{\dimen0=#1 \advance\dimen0 by -#3
       \dimen1=#1 \advance\dimen1 by #3
        \vrule height #1 depth #3 width #3
        \vrule height 0pt depth #3 width #2
        \hskip -#3}}
 \def\RightContraction{\mathord{\, \HollowBox{6pt}{3.1pt}{.4pt}} \kern1.6pt}
\def\qed{{\hfill $\Box$}}
\newtheorem{theorem}{THEOREM}[section]

\newtheorem{lemma}[theorem]{Lemma}

\newtheorem{remark}[theorem]{Remark}
\newtheorem{proposition}[theorem]{Proposition}

\begin{document}
\begin{center}

{\Large \bf Examples of Unbounded Homogeneous \medskip\\
Domains in Complex Space}
 \medskip \\
\normalsize Michael Eastwood\footnote{The first author is supported by the
Australian Research Council.} and Alexander Isaev
\end{center}

\begin{quotation} \small \sl We construct several new examples of homogeneous domains in complex space that do not have bounded realisations. They are equivalent to tubes over affinely homogeneous domains in real space and have a real-analytic everywhere Levi non-degenerate non-umbilic boundary. Using the geometry of the boundary, we determine the full automorphism groups of the domains. We also discuss some interesting examples of tube domains with everywhere umbilic boundary (i.e., boundary equivalent to the corresponding quadric). 
\end{quotation}

\pagestyle{myheadings}
\markboth{M. Eastwood and A. Isaev}{Unbounded Homogeneous Domains}
\setcounter{section}{-1}
\section{Introduction}
\setcounter{equation}{0}

The study of bounded holomorphically homogeneous domains in complex space goes
back to \'E. Cartan \cite{C} who determined all bounded symmetric domains in
$\CC^n$ as well as all bounded homogeneous domains in $\CC^2$ and $\CC^3$. A
fundamental theorem due to Vinberg, Gindikin, and Pyatetskii-Shapiro states
that every bounded homogeneous domain is biholomorphically equivalent to a
Siegel domain of the second kind (see \cite{P-S}). Although this result does
not immediately imply a complete classification of bounded homogeneous domains,
it reduces the classification problem to that for domains of a very special
form. Siegel domains are unbounded by definition (although they possess bounded
realisations) and it is therefore natural to consider not necessarily bounded
homogeneous domains. However, the classification problem in the unbounded case
is extremely hard and is far from complete, despite the existence of a
substantial theory of homogeneous manifolds and spaces (see, e.g., \cite{A}).
Therefore, any new examples of homogeneous domains that possess no bounded
realisations are of interest. In this paper we give such examples.

One way to attempt to produce examples of this sort is to consider tube
domains, that is domains of the form $D_{\Omega}=\Omega+i\RR^n$, where $\Omega$
is a domain in $\RR^n\subset\CC^n$. Such domains are clearly unbounded.
Further, if $\Omega$ is affinely homogeneous, then $D_{\Omega}$ is homogeneous
as a domain in $\CC^n$ since every affine mapping of $\RR^n$ can be lifted to
an affine mapping of $\CC^n$ and since $D_{\Omega}$ is invariant under
imaginary translations. To construct $\Omega$ one can start with an affinely
homogeneous hypersurface $\Gamma\subset\RR^n$ and let $\Omega$ to be a domain
on one side of $\Gamma$. Such a domain $\Omega$ has a chance of being affinely
homogeneous. More precisely, we might hope that the orbits of the affine
symmetry group of $\Gamma$ be, in addition to $\Gamma$ itself, the domains to
either side of it. Of course, in this case, the symmetry group must have
dimension at least $n$ and so its action on $\Gamma$ must have isotropy.
Examples of affinely homogeneous hypersurfaces $\Gamma\subset\RR^n$ can be
taken for instance from the  explicit classifications of affinely homogeneous
curves in $\RR^2$ (see, e.g., \cite{NS2}), surfaces in $\RR^3$ \cite{DKR},
\cite{EE1}, or equiaffinely homogeneous hypersurfaces in $\CC^3$ \cite{NS1} and
$\CC^4$ \cite{EE2}. A classification of affine homogeneous hypersurfaces with
isotropy may be found in~\cite{EE3}. Of course, one must independently verify
that the domain $D_{\Omega}$ so constructed is not biholomorphically equivalent
to any bounded domain and is indeed homogeneous.

The paper is organised as follows. We construct our examples in
Section~\ref{exa}. They are domains in $\CC^4$ arising from domains on each
side of certain affinely homogeneous hypersurfaces
$\Gamma_{\alpha}\subset\RR^4$, where $\alpha$ is a real parameter. We verify
that the domains in question are indeed homogeneous and are not
biholomorphically equivalent to any bounded domain. Further, for a homogeneous
domain it is always desirable to know the full group of holomorphic
automorphisms, and in Section~\ref{autogr} we determine the automorphism groups
of the domains from Section~\ref{exa}. In order to do this, we write the domains in a non-tubular form and study the automorphism group of the boundary which is a real-analytic everywhere Levi non-degenerate non-umbilic hypersurface. In connection with this we recall (see e.g., \cite{R}) that a bounded homogeneous domain with smooth boundary is biholomorphically equivalent to the unit ball. In the unbounded case there are many more domains with smooth boundary, and the examples that we construct in Section~\ref{exa} are not equivalent to any \lq\lq ball-like\rq\rq domains that we discuss further in the paper.

Sometimes a non-trivially looking tube domain in~$\CC^n$, whose boundary is a
tube hypersurface over a homogeneous hypersurface in $\RR^n$, turns out to be
holomorphically equivalent to a simple well-known domain. For example the
domains lying on either side of the quadric
\begin{equation}
\hbox{Re}\,z_n=\langle z',z'\rangle,\label{intro}
\end{equation}
where $\langle z',z'\rangle$ is a Hermitian form in the space of the first
$n-1$ variables $z':=(z_1,\dots,z_{n-1})$ are well-known to be homogeneous. Every such domain possesses a tubular realisation over an affinely homogeneous domain in $\RR^n$. 
Moreover, the convex side of the positive definite quadric admits a bounded
realisation as the the unit ball. For $\alpha=1/12$ the domains that we construct in Section~\ref{exa} are, in fact, of this type. In Section~\ref{quad} we discuss more examples of tube domains
that are holomorphically equivalent to domains lying on one side of
hypersurfaces (\ref{intro}) and announce a theorem that states that the number
of such domains is quite substantial. Therefore, when considering homogeneous
tube domains, one should carefully rule out domains arising in this manner from
quadrics~(\ref{intro}). This is often a non-trivial task. As an example, we
mention a letter by D'Atri (reproduced in the preface to~\cite{G}) where he
apparently found a new homogeneous domain, as one side of the tube over the
Cayley hypersurface~(\ref{cayley}). In fact, as shown in Section~\ref{quad},
the domain that D'Atri considered is equivalent to a domain on one side of the
quadric~(\ref{intro}). In Section~\ref{quad} more examples of this kind are
given.

Before proceeding we would like to acknowledge that this work started while the second author was visiting the University of Adelaide.

\section{The Examples}\label{exa}
\setcounter{equation}{0}

In accordance with the general scheme outlined above we consider a
one-parameter family of affinely homogeneous hypersurfaces in $\RR^4$ that
occurs in the classifications in \cite{EE2} and \cite{EE3}:
\begin{equation}
\Gamma_{\alpha}:=
\left\{(x_1,x_2,x_3,x_4)\in\RR^4:
x_4=x_1x_2+x_3^2+x_1^2x_3+\alpha x_1^4\right\},\alpha\in\RR.\label{gamma}
\end{equation}
We are interested in the domains on each side of $\Gamma_{\alpha}$:
$$
\begin{array}{l}
\Omega_{\alpha}^{>}:=
\left\{(x_1,x_2,x_3,x_4)\in\RR^4:
x_4>x_1x_2+x_3^2+x_1^2x_3+\alpha x_1^4\right\},\\
\vspace{0.2cm}\\
\Omega_{\alpha}^{<}:=
\left\{(x_1,x_2,x_3,x_4)\in\RR^4:
x_4<x_1x_2+x_3^2+x_1^2x_3+\alpha x_1^4\right\}.
\end{array}
$$
One can verify that $\Omega_{\alpha}^{>}$ and $\Omega_{\alpha}^{<}$ are
invariant under the following four subgroups of affine transformations of
$\RR^4$:
$$
\phi_q:\quad
\begin{array}{l}
x_1\mapsto q x_1,\\
x_2\mapsto q^3 x_2,\\
x_3\mapsto q^2 x_3,\\
x_4\mapsto q^4 x_4,
\end{array}
$$
with $q\in\RR^*$,
$$
\psi_r:\quad
\begin{array}{l}
x_1\mapsto x_1+r,\\
x_2\mapsto -4\alpha(4\alpha-1)r^2x_1+x_2+2(4\alpha-1)rx_3
           -\frac{4}{3}\alpha(4\alpha-1)r^3,\\
x_3\mapsto -4\alpha rx_1+x_3-2\alpha r^2,\\
x_4\mapsto -\frac{4}{3}\alpha(4\alpha-1)r^3x_1+rx_2
           +(4\alpha-1)r^2x_3+x_4-\frac{1}{3}\alpha(4\alpha-1)r^4,
\end{array}
$$
with $r\in\RR$,
$$
\mu_s:\quad
\begin{array}{l}
x_1\mapsto x_1,\\
x_2\mapsto x_2+s,\\
x_3\mapsto x_3,\\
x_4\mapsto sx_1+x_4,
\end{array}
$$
with $s\in\RR$,
$$
\nu_t:\quad
\begin{array}{l}
x_1\mapsto x_1,\\
x_2\mapsto -tx_1+x_2,\\
x_3\mapsto x_3+t,\\
x_4\mapsto 2tx_3+x_4+t^2,
\end{array}
$$
with $t\in\RR$.

We will now show that $\Omega_{\alpha}^{>}$ is affinely homogeneous. Take the
point $(0,0,0,1)\in\Omega_{\alpha}^{>}$ and apply the mapping
$F_{q,s,t,r}:=\phi_q\circ\mu_s\circ\nu_t\circ\psi_r$ to it. The result is the
point
$$
\begin{array}{lll}
\Bigl(qr,& q^3(-\frac{4}{3}\alpha(4\alpha-1)r^3-tr+s),& q^2(-2\alpha r^2+t),\\
&&\hspace*{-20pt}
q^4(1-\frac{1}{3}\alpha(4\alpha-1)r^4-4\alpha tr^2+t^2+sr)\Bigr).
\end{array}
$$
Let $(x_1^0,x_2^0,x_3^0,x_4^0)$ be any other point in $\Omega_{\alpha}^{>}$.
Then setting
$$
\begin{array}{l}
\displaystyle q=\left(x_4^0-x_1^0x_2^0-(x_3^0)^2-(x_1^0)^2x_3^0
                        -\alpha (x_1^0)^4\right)^{\frac{1}{4}},\\
\vspace{0.01cm}\\
\displaystyle r=\frac{x_1^0}{q},\\
\vspace{0.01cm}\\
\displaystyle s=\frac{1}{q^3}\left(x_2^0+\frac{4}{3}\alpha(4\alpha-1)(x_1^0)^3
                        +x_1^0x_3^0+2\alpha(x_1^0)^3\right),\\
\vspace{0.01cm}\\
\displaystyle t=\frac{1}{q^2}\left(x_3^0+2\alpha (x_1^0)^2\right),
\end{array}
$$
we obtain an affine automorphism $F_{q,s,t,r}$  of $\Omega_{\alpha}^{>}$ that
maps $(0,0,0,1)$ into $(x_1^0,x_2^0,x_3^0,x_4^0)$. This proves that
$\Omega_{\alpha}^{>}$ is affinely homogeneous. A similar argument shows that
$\Omega_{\alpha}^{<}$ is affinely homogeneous as well.

We will now consider the corresponding tube domains
$D_{\Omega_{\alpha}^{>}},D_{\Omega_{\alpha}^{<}}\subset\CC^4$. Since
$\Omega_{\alpha}^{>}$ and $\Omega_{\alpha}^{<}$ are affinely homogeneous,
$D_{\Omega_{\alpha}^{>}}$ and $D_{\Omega_{\alpha}^{<}}$ are holomorphically
homogeneous. We may write these domains in a different form. Suppose firstly
that $\alpha\ne 1/12$. Then the mapping
$$
\begin{array}{l}
\displaystyle z_1\mapsto \left|\frac{3}{2}\left(\alpha
                         -\frac{1}{12}\right)\right|^{\frac{1}{4}}z_1,\\
\vspace{0.001cm}\\
\displaystyle z_2\mapsto \frac{z_2+z_1z_3+\alpha z_1^3}
 {\left|\frac{3}{2}\left(\alpha-\frac{1}{12}\right)\right|^{\frac{1}{4}}},\\
\vspace{0.001cm}\\
\displaystyle z_3\mapsto \sqrt{2}\left(z_3+\frac{z_1^2}{4}\right),\\
\vspace{0.001cm}\\
\displaystyle z_4\mapsto 4z_4-2z_1z_2-2z_3^2-z_1^2z_3-\frac{\alpha}{2}z_1^4
\end{array}
$$
transforms $D_{\Omega_{\alpha}^{>}}$ and $D_{\Omega_{\alpha}^{<}}$ into
\begin{equation}
D^{>}_{+}:=
\left\{(z_1,z_2,z_3,z_4)\in\CC^4:
\hbox{Re}\,z_4>z_1\overline{z_2}+z_2\overline{z_1}+|z_3|^2+|z_1|^4\right\}
\label{gp}
\end{equation}
and
\begin{equation}
D^{<}_{+}:=
\left\{(z_1,z_2,z_3,z_4)\in\CC^4:
\hbox{Re}\,z_4<z_1\overline{z_2}+z_2\overline{z_1}+|z_3|^2+|z_1|^4\right\}
\label{lp}
\end{equation}
respectively, if $\alpha>1/12$, and into
\begin{equation}
D^{>}_{-}:=
\left\{(z_1,z_2,z_3,z_4)\in\CC^4:
\hbox{Re}\,z_4>z_1\overline{z_2}+z_2\overline{z_1}+|z_3|^2-|z_1|^4\right\}
\label{gm}
\end{equation}
and
\begin{equation}
D^{<}_{-}:=
\left\{(z_1,z_2,z_3,z_4)\in\CC^4:
\hbox{Re}\,z_4<z_1\overline{z_2}+z_2\overline{z_1}+|z_3|^2-|z_1|^4\right\}
\label{lm}
\end{equation}
respectively, if $\alpha<1/12$. If, however, $\alpha=1/12$, then the mapping
 \begin{equation}
\begin{array}{l}
\displaystyle z_1\mapsto \frac{1}{\sqrt{2}}
                   \left(z_1+z_2+z_1z_3+\frac{1}{12}z_1^3\right),\\
\vspace{0.001cm}\\
\displaystyle z_2\mapsto \sqrt{2}\left(z_3+\frac{z_1^2}{4}\right),\\
\vspace{0.001cm}\\
\displaystyle z_3\mapsto  \frac{1}{\sqrt{2}}
                    \left(z_1-z_2-z_1z_3-\frac{1}{12}z_1^3\right),\\
\vspace{0.001cm}\\
\displaystyle z_4\mapsto 4z_4-2z_1z_2-2z_3^2-z_1^2z_3-\frac{1}{24}z_1^4
\end{array}\label{quadequiv}
\end{equation}
transforms $D_{\Omega_{1/12}^{>}}$ and $D_{\Omega_{1/12}^{<}}$ into
\begin{equation}
D^{>}_{0}:=
\left\{(z_1,z_2,z_3,z_4)\in\CC^4:
\hbox{Re}\,z_4>|z_1|^2+|z_2|^2-|z_3|^2\right\}\label{g0}
\end{equation}
and
\begin{equation}
D^{<}_{0}:=
\left\{(z_1,z_2,z_3,z_4)\in\CC^4:
\hbox{Re}\,z_4<|z_1|^2+|z_2|^2-|z_3|^2\right\}\label{l0}
\end{equation}
respectively.

We discuss the domains $D^{>}_0$ and $D^{<}_0$ in greater generality in Section
\ref{quad} and for the moment concentrate on the domains $D^{>}_{\pm}$,
$D^{<}_{\pm}$. We note that none of these domains is biholomorphically
equivalent to a bounded domain. In fact, all these domains are not Kobayashi-hyperbolic (we remark here that it is shown in \cite{N} that any connected homogeneous Kobayashi-hyperbolic manifold is biholomorphically equivalent to a bounded domain in complex space). Indeed, domains $D^{>}_{+}$, $D^{>}_{-}$ contain the
affine complex line $\{z_1=0,z_3=0,z_4=1\}$ and the domains $D^{<}_{+}$,
$D^{<}_{-}$ the affine complex line $\{z_1=0,z_3=0,z_4=-1\}$. Hence, we have
proved the following

\begin{theorem}\label{main}\sl
The domains $D^{>}_{\pm}, D^{<}_{\pm}\subset\CC^4$ defined by
(\ref{gp})--(\ref{lm}) are holomorphically homogeneous and not
biholomorphically equivalent to any bounded domain.
\end{theorem}

\section{The Automorphism Groups of $D^{>}_{\pm}$
and $D^{<}_{\pm}$}\label{autogr}
\setcounter{equation}{0}

In this section we determine the groups $\hbox{Aut}\,(D^{>}_{\pm})$ and $\hbox{Aut}\,(D^{<}_{\pm})$ of holomorphic automorphisms of $D^{>}_{\pm}$ and $D^{<}_{\pm}$ respectively.
Let $P_{\pm}$ be the following subgroup of holomorphic transformations of $\CC^4$:
\begin{equation}
\begin{array}{lll}
z_1&\mapsto &qe^{i\phi}z_1+\rho,\\
\vspace{0.1cm} &&\\
z_2&\mapsto &\displaystyle
(\mp2|\rho|^2qe^{i\phi}+q^2b)z_1+q^3e^{i\phi}z_2+qdz_3
\mp2\overline{\rho}q^2e^{2i\phi}z_1^2+\sigma,\\
\vspace{0.1cm} &&\\
z_3&\mapsto &\displaystyle -\overline{d}e^{i(\phi+\psi)}z_1
+q^2e^{i\psi}z_3+\tau,\\
\vspace{0.1cm} &&\\
z_4&\mapsto &(2\overline{\sigma}qe^{i\phi}+2\overline{\rho}q^2b
-2\overline{\tau}\overline{d}e^{i(\phi+\psi)})z_1
+2\overline{\rho}q^3e^{i\phi}z_2+{}\\
\vspace{0.1cm} &&\\
&&(2\overline{\rho}qd+
2\overline{\tau}q^2e^{\i\psi})z_3
+q^4z_4\mp2\overline{\rho}^2q^2e^{2i\phi}z_1^2+\rho\overline{\sigma}
+\sigma\overline{\rho}+{}\\
\vspace{0.1cm} &&\\
&&|\tau^2|\pm|\rho^4|+iu,
\end{array}\label{thegroup}
\end{equation}
where $q>0$, $\phi,\psi,u\in\RR$, $\rho,\sigma,\tau,b,d\in\CC$,
$\hbox{Re}(e^{i\phi}\overline{b})\le 0$ and
\begin{equation}
|d|^2=-2q^3\hbox{Re}(e^{i\phi}\overline{b}).\label{thecondition}
\end{equation}
It can be checked directly that $P_{\pm}$ preserves $D_{\pm}^{>}$ and
$D_{\pm}^{<}$, so $P_{\pm}\subset \hbox{Aut}(D_{\pm}^{>})$ and
$P_{\pm}\subset \hbox{Aut}(D_{\pm}^{<})$.

Below we  prove the following theorem.

\begin{theorem}\label{domgroup}\sl
$\hbox{Aut}(D_{\pm}^{>})=\hbox{Aut}(D_{\pm}^{<})=P_{\pm}$.
\end{theorem}

To prove Theorem \ref{domgroup} we deal with the automorphism group of $M_{\pm}:=\partial D_{\pm}^{>}=\partial D_{\pm}^{<}$. Denote by $\hbox{Aut}(M_{\pm})$ the group of CR-automorphisms of $M_{\pm}$ equipped with the topology of uniform convergence of the derivatives of all orders of the component functions on compact subsets of $M_{\pm}$. The Levi form of $M_{\pm}$ at every point is non-degenerate, and therefore $\hbox{Aut}(M_{\pm})$ is a Lie group in this topology \cite{T}. We need the following theorem that implies Theorem \ref{domgroup}. 

\begin{theorem}\label{w}\sl $\hbox{Aut}(M_{\pm})=P_{\pm}$.
\end{theorem}

We first show how Theorem \ref{domgroup} follows from Theorem \ref{w}.
Since the Levi form of $M_{\pm}$ at every point has eigenvalues of opposite
signs, every element $f$ of $\hbox{Aut}(D_{\pm}^{>})$ or
$\hbox{Aut}(D_{\pm}^{<})$ extends past $M_{\pm}$ to a biholomorphic map between
neighbourhoods of $\overline{D_{\pm}^{>}}$ and $\overline{D_{\pm}^{<}}$
respectively thus giving rise to an element  $g\in\hbox{Aut}(M_{\pm})$.
Clearly, $f$ is uniquely determined by $g$. Since
$\hbox{Aut}(M_{\pm})=P_{\pm}$ by Theorem \ref{w}, we obtain
$\hbox{Aut}(D_{\pm}^{>})=\hbox{Aut}(D_{\pm}^{<})=P_{\pm}$ as required.
\smallskip\\

\noindent {\bf Proof of Theorem \ref{w}:} First we note that $P_{\pm}$
considered as a subgroup of $\hbox{Aut}(M_{\pm})$ is closed and is therefore a
Lie subgroup of $\hbox{Aut}(M_{\pm})$. Clearly, the topology induced on
$P_{\pm}$ by $\hbox{Aut}(M_{\pm})$ coincides with that induced on $P_{\pm}$ by
its parameters as in formula (\ref{thegroup}), and hence
the dimension of $P_{\pm}$ as a subgroup of $\hbox{Aut}(M_{\pm})$ is equal to 13. Therefore $\hbox{dim}\,\hbox{Aut}(M_{\pm})\ge 13$.    

We will now show that $\hbox{dim}\,\hbox{Aut}(M_{\pm})\le 13$. Since $M_{\pm}$
is homogeneous, it is sufficient to prove that $\hbox{dim}\,I_0(M_{\pm})\le 6$,
where $I_0(M_{\pm})\subset \hbox{Aut}(M_{\pm})$ is the isotropy subgroup of the
point $0\in M_{\pm}$. More generally, we prove the following proposition.

\begin{proposition}\label{dim}\sl Let $S\subset\CC^4$ be a real-analytic hypersurface passing through the origin and assume that the Levi form of $S$ is everywhere non-degenerate and has signature (2,1). Suppose further that $S$ is non-umbilic at 0. Let $I_0(S)$ denote the Lie group of all CR-automorphisms of $S$ preserving 0. Then $\hbox{dim}\,I_0(S)\le 6$.
\end{proposition}

\noindent {\bf Proof of Proposition \ref{dim}:} We will use the Chern-Moser
normal form \cite{CM} for $S$. Namely, we will say that $S$ in coordinates
$z=(z_1,z_2,z_3),w=u+iv$ near $0$ is written in the Chern-Moser normal form if
$S$ near $0$ is given by an equation
\begin{equation}
u=\langle z,z\rangle
+\sum_{k,\overline{l}\ge 2}F_{k\overline{l}}(z,\overline{z},v),
\label{normal}
\end{equation}
where $\langle z,z\rangle$ is a non-degenerate Hermitian form of signature
$(2,1)$ and\linebreak $F_{k\overline{l}}(z,\overline{z},v)$ are polynomials of
degree $k$ in $z$ and $\overline{l}$ in $\overline{z}$  whose coefficients are
analytic functions of $v$ such that the following conditions hold
$$
\begin{array}{l}
\hbox{tr}\,F_{2\overline{2}}\equiv0,\\
\hbox{tr}^2\,F_{2\overline{3}}\equiv0,\\
\hbox{tr}^3\,F_{3\overline{3}}\equiv0,
\end{array}
$$
with the operator $\hbox{tr}$ defined as follows. Let $\langle z,z\rangle=\sum_{\alpha,\beta=1}^4h^{\alpha\beta}z_{\alpha}\overline{z_{\beta}}$ and let $(g_{\alpha\beta})$ be the matrix inverse to $(h^{\alpha\beta})$. Then we have
$$
\displaystyle \hbox{tr}:=\sum_{\alpha,\beta=1}^4g_{\alpha\beta}\frac{\partial^2 }{\partial z_{\alpha}\partial \overline{z_{\beta}}}.
$$

Let $J_0(S)$ denote the group of all local CR-automorphisms of $S$ defined near
the origin and preserving it. It follows from \cite{E} that near $0$ one can
choose holomorphic coordinates in which $S$ is given in the Chern-Moser normal
form and every element of $J_0(S)$ is written as
\begin{equation}
\begin{array}{l}
z\mapsto \lambda Uz,\\
w\mapsto \lambda^2 w,
\end{array}\label{linear}
\end{equation}
where $U$ is a matrix such that $\langle Uz,Uz\rangle=\langle z,z\rangle$ and
$\lambda>0$. Without loss of generality we assume that $\langle z,z\rangle=|z_1|^2+|z_2|^2-|z_3|^2$
and hence the matrix $U$ satisfies
\begin{equation}
U^tH\overline{U}=H,\label{group}
\end{equation}
where
$$
H:=\left(
\begin{array}{rrr}
1 & 0 & 0\\
0 & 1 & 0\\
0 & 0 & -1
\end{array}
\right).
$$
The group $G$ of all matrices given by condition (\ref{group}) is isomorphic to
the usual group $U(2,1)$ of pseudo-unitary matrices $V$ satisfying
$$
VHV^*=H
$$
by means of the mapping
$$
V=(U^t)^{-1}.
$$

It follows from \cite{B}, \cite{L} that, for every element of $J_0(S)$, $\lambda$
is uniquely determined by the corresponding matrix $U$ by means of an algebraic
relation. This implies that the subgroup $G_0\subset G$ of matrices $U$ arising
from automorphisms in $J_0(S)$ is a closed subgroup of $G$. In addition, the mapping $U\mapsto\lambda$
is a Lie group homomorphism from $G_0$ into~$\RR^*$.

We will now plug an automorphism of the form (\ref{linear}) into equation
(\ref{normal}). We obtain:
\begin{equation}
\displaystyle \frac{1}{\lambda^2}
\sum_{k,\overline{l}\ge 2}
F_{k\overline{l}}(\lambda Uz,\overline{\lambda Uz},\lambda^2v)=
\sum_{k,\overline{l}\ge 2}F_{k\overline{l}}(z,\overline{z},v).\label{raw1}
\end{equation}
Since $S$ is not umbilic at 0, $F_{2\overline{2}}(z,\overline{z},0)\not\equiv
0$. Extracting from (\ref{raw1}) terms independent of $v$ of degree 2 in each
of $z$ and $\overline{z}$ we obtain:
\begin{equation}
\displaystyle F_{2\overline{2}}(Uz,\overline{Uz},0)=
\frac{1}{\lambda^2} F_{2\overline{2}}(z,\overline{z},0).\label{final}
\end{equation}

Hence a necessary condition for a matrix $U\in G$ to belong to $G_0$ is the
preservation of the term $F_{2\overline{2}}(z,\overline{z},0)$ up to a scalar
multiple as in (\ref{final}). We will show that condition (\ref{final}) implies
that $\hbox{dim}\,G_0\le 6$.

Suppose first that $\hbox{dim}\,G_0=9$, i.e., $G_0=G$. Since there does not
exist a non-trivial homomorphism from $G$ into $\RR^*$, we have $\lambda=1$ for
all $U$. Then (\ref{final}) gives that $F_{2\overline{2}}(z,\overline{z},0)$ is
a function of $\langle z,z\rangle$, i.e.\
$F_{2\overline{2}}(z,\overline{z},0)=c\langle z,z\rangle^2$, for some
$c\in\RR^*$. This
is impossible since then
$\hbox{tr}\,F_{2\overline{2}}(z,\overline{z},v)
=8c\langle z,z\rangle\not\equiv 0$.

To deal with the cases $\hbox{dim}\,G_0=7,8$ we need the following lemma.

\begin{lemma}\label{subgroups}\sl\mbox{ }
\begin{flushleft}\begin{tabular}{rl}
(i) & The only closed subgroup of $U(2,1)$ is codimension 1 is $SU(2,1)$.\\
(ii)&There does not exist a closed subgroup of $U(2,1)$ of codimension 2.
\end{tabular}\end{flushleft}
\end{lemma}

\noindent {\bf Proof of Lemma \ref{subgroups}:} It suffices to prove this on
the level of Lie algebras, which we shall denote by ${\frak{u}}(2,1)$ and
${\frak{su}}(2,1)$, respectively. Certainly, ${\frak{su}}(2,1)$, has
codimension 1 in ${\frak{u}}(2,1)$ and any other subalgebra of codimension 1 or
2 would intersect ${\frak{su}}(2,1)$ in a subalgebra of codimension 1 or 2.
Therefore, it suffices to show that ${\frak{su}}(2,1)$ has no subalgebras of
codimension 1 or~2, equivalently of dimension 7 or~6. To do this, we shall show
that the complex Lie algebra ${\frak{sl}}(3,\CC)$ has no complex subalgebras of
dimension 7 and classify those of dimension~6. The result concerning
${\frak{su}}(2,1)$ will follow if we can show that these 6-dimensional
complex subalgebras have no real form in ${\frak{su}}(2,1)$.

Recall that the Killing form
$$\langle X,Y\rangle={\rm{trace}}\,(XY),\quad\mbox{for }X,Y\in
{\frak{sl}}(3,{\CC}),$$
is a non-degenerate symmetric form. Therefore, if
${\frak{s}}\subset{\frak{sl}}(3,{\CC})$ is a subalgebra, the linear subspace
${\frak{s}}^\perp$ with respect to the Killing form, will have dimension equal
to the codimension of~${\frak{s}}$. Furthermore, invariance of the Killing form
with respect to the adjoint representation implies that, if $P\in{\frak{s}}^\perp$,
then the linear mapping
\begin{equation}\label{strong}
[P,{}\cdot{}]:{\frak{s}}\to{\frak{sl}}(3,\CC)
\end{equation}
has range contained in~${\frak{s}}^\perp$. Now, if
$\dim{\frak{s}}^\perp\leq 2$, then the dimension of the kernel of
$[P,{}\cdot{}]$ must be at least 4 and this will prove to be
very restrictive. In particular, since ${\frak{s}}\subset P^\perp$, we conclude
that
\begin{equation}\label{weak}
[P,{}\cdot{}]:P^\perp\to{\frak{sl}}(3,\CC)\quad
\mbox{has kernel of dimension }\geq 4.
\end{equation}
Notice that this constraint depends only on the element
$P\in{\frak{sl}}(3,\CC)$ and so we may test it for any particular~$P$.
Furthermore, we may assume without loss of generality that $P$ is in Jordan
canonical form in which case a simple computation shows that there is only one
$P$ satisfying~(\ref{weak}), namely
$$P=\left(\begin{array}{ccc}0&1&0\\ 0&0&0\\ 0&0&0\end{array}\right).$$
This immediately rules out subalgebras of codimension 1 since $P^\perp$ is not
a subalgebra. 

The constraint on mapping (\ref{strong}) now pins down ${\frak{s}}$ as two
possibilities, namely matrices of the form either
$$\left(\begin{array}{ccc}*&*&*\\ 0&*&*\\ 0&*&*\end{array}\right)\quad\mbox{or}
\quad\left(\begin{array}{ccc}*&*&*\\ 0&*&0\\ {*}&*&*\end{array}\right).
$$
Both of these are, indeed, subalgebras and have geometric interpretations. The
first is the stabiliser up to scale of the first standard basis vector in the
defining representation. If it were the complexification of a subalgebra
$${\frak{s}}_0\subset{\frak{su}}(2,1)\subset{\frak{sl}}(3,\CC)$$
then ${\frak{s}}_0$ would stabilise up to scale some vector~$v\in\CC$. Up to
conjugation and scale, there are only three possibilities for $v$ according to
$\langle v,v\rangle>0$, $\langle v,v\rangle<0$, or $\langle v,v\rangle=0$. Taking $v$ in some convenient normal
form it is easy to check that the corresponding stabiliser has dimension only
$4$, $4$, or $5$ respectively. The second possibility for
${\frak{s}}\subset{\frak{sl}}(3,\CC)$ is similarly eliminated thanks to a
geometric interpretation in the dual representation.

The lemma is proved.\qed

\begin{remark}\label{all}\rm This proof of Lemma~\ref{subgroups} extends to
higher dimensions, where it yields the subalgebras of ${\frak{sl}}(n,\CC)$ of
maximal dimension. They are parabolic and their real forms may be determined
following the Satake classification. More generally, the maximal subalgebras of
the complex classical Lie algebras were determined by Dynkin~\cite{D}. The real
case was considered by Komrakov~\cite{K} and, in principle, the lemma follows
from his classification. However, no proofs are given in~\cite{K}. One may
construct another proof by considering how potential subalgebras of
${\frak{su}}(2,1)$ intersect ${\frak{u}}(2)$, whereupon Lemma~2.1 of~\cite{IK},
dealing with large compact subgroups of $GL(n,\CC)$, may be used to eliminate
the various possibilities. Finally, it was pointed out to us by Vladimir Ezhov
that a careful investigation of the proof of his linearisation result
in~\cite{E}, shows that, for a hypersurface of the form
$$
\hbox{Re}\,z_4=z_1\overline{z_2}+z_2\overline{z_1}+|z_3|^2\pm |z_1|^4
$$
as we have, any local CR-automorphism near the origin is already linear. By
this means one can avoid Lemma~\ref{subgroups} if one so chooses.
\end{remark}

We now finish the proof of Proposition \ref{dim}. Suppose that $\hbox{dim}\,G_0=8$. Then by Lemma \ref{subgroups},  $G_0$ is the subgroup of $G$ given by the condition $\hbox{det}\,U=1$. In this case we again get that     $F_{2\overline{2}}(z,\overline{z},0)$ is a function of $\langle z,z\rangle$ and obtain a contradiction as before. Further, Lemma \ref{subgroups} gives that $\hbox{dim}\,G_0\ne7$, and hence $\hbox{dim}\,G_0\le 6$.

Finally, since the mapping $f\mapsto U$ is an injective Lie group homomorphism from $I_0(S)$ into $G_0$, we obtain that $\hbox{dim}\,I_0(S)\le 6$, and the proposition is proved.\qed 
\begin{remark}\rm Instead of the non-umbilicity of $S$ at 0 it is sufficient to
assume in Proposition \ref{dim} that at least one of $F_{2\overline{2}}$,
$F_{2\overline{3}}$ and $F_{3\overline{3}}$ is not identically zero.
\end{remark}

We will now continue with the proof of Theorem \ref{w}. Proposition \ref{dim}
gives that $\hbox{dim}\,\hbox{Aut}(M_{\pm})\le 13$ and hence in fact
$\hbox{dim}\,\hbox{Aut}(M_{\pm})=13$. Since $P_{\pm}$ is connected in
$\hbox{Aut}(M_{\pm})$, it coincides with the connected component of the
identity of $\hbox{Aut}(M_{\pm})$. To finish the proof of Theorem \ref{w} we
need the following proposition.

\begin{proposition}\label{connected}\sl The group $\hbox{Aut}(M_{\pm})$ is
connected.
\end{proposition}

\noindent{\bf Proof of Proposition \ref{connected}:} Since $M_{\pm}$ is
homogeneous and connected, it is sufficient to prove that $I_0(M_{\pm})$ is
connected. As in the proof of Proposition \ref{dim}, let $J_0(M_{\pm})$ denote the group
of all local CR-automorphisms of $M_{\pm}$ defined near the origin and
preserving it. Let $f\in J_0(M_{\pm})$. Since $f$ extends holomorphically to a
neighbourhood of the origin, we can write $f$ in the form:
$$
\begin{array}{l}
z^*=f_1(z,w),\\
w^*=f_2(z,w),
\end{array}
$$
where $f_1,f_2$ are holomorphic in a neighbourhood of the origin in $\CC^4$.
Let
$$
\displaystyle U_f:=
\frac{1}{\sqrt{\partial f_2/\partial w (0)}}\frac{\partial f_1}{\partial z}(0),
$$
(note that $\partial f_2/\partial w (0)$ is necessarily positive). Clearly,
$U_f$ is a matrix satisfying
\begin{equation}
\langle U_fz,U_fz\rangle=\langle z,z\rangle,\label{unit}
\end{equation}
where $\langle z,z\rangle:=z_1\overline{z_2}+z_2\overline{z_1}+|z_3|^2$. It
follows from \cite{B}, \cite{L} that $f$ is uniquely determined by the
corresponding $U_f$ and that the collection of all matrices $U_f$ arising in
this way from elements of $J_0(M_{\pm})$ form a closed subgroup $G_0$ in the group $G$
of all matrices satisfying (\ref{unit}) (of course, $G$ is isomorphic to
$U(2,1)$).

Since $M_{\pm}$ is non-umbilic at 0, it follows from the proof of
Proposition~\ref{dim} that $\hbox{dim}\,G_0\le 6$. Let $G_0'$ be the subgroup
of $G_0$ that consists of matrices $U_f$ for $f\in I_0(M_{\pm})$. We will prove
the connectedness of $I_0(M_{\pm})$ by showing that $G_0'$ is connected. It
follows from (\ref{thegroup}) that $G_0'$ contains matrices of the following
form:
\begin{equation}
\left(
\begin{array}{ccc}
\displaystyle\frac{1}{q} e^{i\phi} & 0 & 0\\
b & qe^{i\phi} & \displaystyle\frac{d}{q}\\
-\displaystyle\frac{\overline{d}}{q^2}e^{i(\phi+\psi)} & 0 & e^{i\psi}
\end{array}
\right),\label{linn}
\end{equation}
where $q>0$, $\phi,\psi\in\RR$, $b,d\in\CC$,
$\hbox{Re}(e^{i\phi}\overline{b})\le 0$ and condition (\ref{thecondition})
holds. The group $G_0''$ of matrices (\ref{linn}) is a closed connected
subgroup of $G$ of dimension 6. Hence $G_0''$ is the connected component of the
identity of $G_0'$.

We now need the following lemma.

\begin{lemma}\label{subb}\sl There does not exist a disconnected (not necessarily closed) subgroup of $G$ whose connected component of the identity coincides with $G_0''$.
\end{lemma}

\noindent {\bf Proof of Lemma \ref{subb}:} Let $H\subset G$ be a group whose connected component of the identity is $G_0''$. Then any other connected component of $H$ is of the form $gG_0''$ for some $g\in G$, so we have $H=\cup_{\alpha} g_{\alpha}G_0''$. For every pair of indices $\alpha,\beta$ there exists an index $\gamma$ such that
$$
g_{\alpha}G_0''g_{\beta}G_0''=g_{\gamma}G_0'',
$$
or
$$
g_{\alpha}G_0''g_{\beta}=g_{\gamma}G_0''.
$$
Choosing $g_{\alpha}G_0''=G_0''$ we obtain
\begin{equation}
G_0''g_{\beta}=g_{\gamma}G_0''.\label{comp}
\end{equation}
We now apply both sides of (\ref{comp}) to the vector $v:=(0,1,0)$. For every $g\in G_0''$ we have $gv=\lambda(g)v$ with $\lambda(g)\in\CC^*$, and therefore for every $g\in G_0''$ we obtain
$$
g(g_{\beta}v)=\lambda(g)(g_{\gamma}v),
$$
or, denoting $w_1:=g_{\alpha}v$ and $w_2:=g_{\gamma}v$,
$$
gw_1=\lambda(g) w_2,
$$
i.e., the whole group $G_0''$ maps the vector $w_1$ into the complex line generated by $w_2$. It is then easy to see that such a vector $w_1$ has to be proportional to $v$. 

Hence we obtain that $g_{\beta}$ preserves $v$ up to a multiple. This implies that $g_{\beta}$ has the form (\ref{linn}), i.e.,  $g_{\beta}\in G_0''$ for all $\beta$, and $H=G_0''$.

The lemma is proved.\qed

We can now finish the proof of Proposition \ref{connected}. Indeed, Lemma \ref{subb} implies that $G_0'=G_0''$ is connected. Since the mapping $f\mapsto U_f$ is an injective Lie group homomorphism from $I_0(M_{\pm})$ onto $G_0'$ and since the image of the connected component of the identity of $I_0(M_{\pm})$ under this mapping is $G_0'$, we obtain that $I_0(M_{\pm})$ is connected.\qed

Further, since $\hbox{Aut}(M_{\pm})$ is connected by Proposition \ref{connected}, we obtain that $\hbox{Aut}(M_{\pm})=P_{\pm}$, and Theorem \ref{w} is proved.\qed  

\section{Domains with Boundary Equivalent to the Quadric}\label{quad}
\setcounter{equation}{0}

We will now turn to domains $D^{>}_0$ and $D^{<}_0$ defined in (\ref{g0}) and
(\ref{l0}). In fact, they belong to the following well-known class of domains
in $\CC^n$. Let $(z_1,\dots,z_n,z_{n+1})$ be coordinates in $\CC^{n+1}$, and
$x_j=\hbox{Re}\,z_j$, $j=1,\dots,n+1$.  Set $z:=(z_1,\dots,z_n)$ and
$x:=(x_1,\dots,x_n)$. Consider a non-degenerate Hermitian form
$H_{p,n}(z,\overline{z})$ on $\CC^n$, where $p$ is the number of positive
eigenvalues of $H_{p,n}$ and suppose that $n\le 2p$. Without loss of generality
we assume that $H_{p,n}$ is given in the diagonal form
$$
H_{p,n}(z,\overline{z})=\sum_{j=1}^p|z_j|^2-\sum_{j=p+1}^n|z_j|^2.
$$
We now set
$$
\begin{array}{l}
D_{H_{p,n}}^{>}:=
\{(z,z_{n+1})\in\CC^{n+1}:x_{n+1}>H_{p,n}(z,\overline{z})\},\\
D_{H_{p,n}}^{<}:=\{(z,z_{n+1})\in\CC^{n+1}:x_{n+1}<H_{p,n}(z,\overline{z})\}.
\end{array}
$$
It is easy to check that the mappings
$$
\begin{array}{l}
z\mapsto az+b,\\
z_{n+1}\mapsto 2aH_{p,n}(z,\overline{b})+a^2z_{n+1}+H_{p,n}(b,\overline{b})+ic,
\end{array}
$$
with $a\in\RR^*$, $b\in\CC^n$, $c\in \RR$, act transitively on each of
$D_{H_{p,n}}^{>}$ and $D_{H_{p,n}}^{<}$ so these domains are homogeneous. The
domain $D_{H_{p,n}}^{<}$ is never Kobayashi-hyperbolic as it contains the
complex line $\{z_2=\dots=z_n=0,z_{n+1}=-1\}$.  The domain $D_{H_{p,n}}^{>}$ is
not Kobayashi hyperbolic for $p<n$ as it contains the
complex line $\{z_1=\dots=z_{n-1}=0,z_{n+1}=1\}$. The domain $D_{H_{n,n}}^{>}$
is biholomorphically equivalent to the unit ball.

Thus the domains $D_{\Omega_{1/12}^{>}}$ and $D_{\Omega_{1/12}^{<}}$ serve as a
warning that sometimes non-trivially looking homogeneous tube domains can in
fact be biholomorphically equivalent to a well-known domain (in our example
$D^{>}_0$ and $D^{<}_0$).

It is interesting to remark that $D_{H_{p,n}}^{>}$ and $D_{H_{p,n}}^{<}$ can
always be realised as tube domains. The mapping
$$
\begin{array}{l}
z\mapsto\sqrt{2}z,\\
z_{n+1}\mapsto z_{n+1}+H_{p,n}(z,z)
\end{array}
$$
transforms $D_{H_{p,n}}^{>}$ and $D_{H_{p,n}}^{<}$ into the tube domains
$$
\left\{(z,z_{n+1})\in\CC^{n+1}:x_{n+1}>H_{p,n}(x,x)\right\}
$$
and
$$
\left\{(z,z_{n+1})\in\CC^{n+1}:x_{n+1}<H_{p,n}(x,x)\right\}
$$
respectively. The bases of these domains are affinely homogeneous in
$\RR^{n+1}$. In \cite{DY}, \cite{I1}, \cite{IM} all tube hypersurfaces locally
equivalent to the quadric
$$
Q_{H_{p,n}}:=\left\{(z,z_{n+1})\in\CC^{n+1}:x_{n+1}=H_{p,n}(z,z)\right\}
$$
were determined for $n-p\le 2$. Among such hypersurfaces those given by
polynomial graphs are always globally equivalent to $Q_{H_{p,n}}$, and their
bases are affinely homogeneous in $\RR^{n+1}$. Moreover, the domains on each
side of such a polynomial graph are affinely homogeneous.

For $p=n$ the only polynomial tube hypersurface (up to affine equivalence) is
$$
\left\{(z,z_{n+1})\in\CC^{n+1}:x_{n+1}=H_{n,n}(x,x)\right\}.
$$

As $n-p$ grows, higher order polynomial hypersurfaces appear.  The tube over
the hypersurface $\Gamma_{1/12}$ defined in (\ref{gamma}) is an example for
$n=3, p=2$, and mapping (\ref{quadequiv}) establishes equivalence between this
tube and $Q_{H_{2,3}}$  (see also \cite{IM}). Another non-trivial example
occurs for $n=2, p=1$. Consider the hypersurface in $\RR^3$ given by the
equation:
\begin{equation}
x_3=x_1x_2+x_1^3.\label{cayley}
\end{equation}
It is called the Cayley surface (see \cite{NS2} for a discussion of its
properties). The Cayley surface is affinely homogeneous, and the domains on
each side of it are affinely homogeneous as well. The tube over the Cayley
surface appears in the classification in \cite{IM}, where its equivalence to
$Q_{H_{1,2}}$ was explicitly established. Namely, the mapping
$$
\begin{array}{l}
\displaystyle z_1\mapsto
\frac{1}{\sqrt{2}}\left(z_1+z_2+\frac{3}{2}z_1^2\right),\\
\vspace{0.001cm}\\
\displaystyle z_2\mapsto
\frac{1}{\sqrt{2}}\left(z_1-z_2-\frac{3}{2}z_1^2\right),\\
\vspace{0.001cm}\\
z_3\mapsto 4z_3-2z_1z_2-z_1^3
\end{array}
$$
transforms the tube over the Cayley surface into
$$
\left\{(z_1,z_2,z_3)\in\CC^3:\hbox{Re}\,z_3=|z_1|^2-|z_2|^2\right\}.
$$
For $n-p=2$ very complicated polynomial hypersurfaces equivalent to
$Q_{H_{p,n}}$ occur. For example, if $n=7, p=5$ we have the following
one-parameter family of pairwise affinely non-equivalent hypersurfaces
\cite{I1}:
$$
\begin{array}{l}
\displaystyle x_8
=x_1^2+x_2^2+x_3^2+x_4x_5+x_6x_7+{}\\
\phantom{x_8={}}\displaystyle 2\sqrt{2(1+\sigma)}x_1x_4x_6
+2\sqrt{3\sigma}x_2x_6^2+\frac{1+\sigma}{\sqrt{3\sigma}}x_2x_4^2+{}\\
\phantom{x_8={}}\displaystyle\sqrt{\frac{-\sigma^2
+34\sigma-1}{3\sigma}}x_3x_4^2+(x_4^2+x_6^2)(x_4^2+\sigma x_6^2),\quad
\sigma\in[1,17+12\sqrt{2}).
\end{array}
$$
Every such a hypersurface and the domains on each side of it are affinely
homogeneous, and it is possible to write explicitly a polynomial biholomorphism
that maps the tubes over these hypersurfaces onto $Q_{H_{5,7}}$.  In general,
the following holds.

\begin{theorem}\label{spher}\sl For every $p\le n$ and every $2\le s\le 2(n-p)+2$ there exists a polynomial $P(x_1,\dots,x_n)$ of degree $s$ such that:
\smallskip\\

\noindent (i) the graph of $P$ is an affinely homogeneous hypersurface in $\RR^{n+1}$;
\smallskip\\

\noindent (ii) the tube hypersurface over the graph of $P$ is equivalent to $Q_{H_{p,n}}$ by means of a polynomial mapping.
\smallskip\\

Moreover, if $n\ge 7$ and $n-p\ge 2$, then for some $2\le s_0\le 2(n-p)+2$ there exists a family $\{P_{\sigma}\}$ of polynomials of degree $s_0$ depending on a continuous parameter, that possess properties (i) and  (ii), such that the graphs of $P_{\sigma}$ are pairwise affinely non-equivalent.
\end{theorem}

We do not prove Theorem \ref{spher} here; it follows from the techniques
developed in \cite{I2}. Theorem \ref{spher} shows that one has to exercise
caution when constructing homogeneous domains from tubes over affinely
homogeneous hypersurfaces in real space since a substantial number of such
domains are either $D_{H_{p,n}}^{>}$, or $D_{H_{p,n}}^{<}$ in disguise.

{\obeylines
Department of Pure Mathematics
University of Adelaide
South Australia 5005
AUSTRALIA
E-mail: meastwoo@maths.adelaide.edu.au
\hbox{ \ \ }

\hbox{ \ \ }
Department of Mathematics
The Australian National University
Canberra, ACT 0200
AUSTRALIA
E-mail: alexander.isaev@maths.anu.edu.au}


\begin{thebibliography}{ABCD}

\bibitem[A]{A} Akhiezer, D. N., Homogeneous complex manifolds
  (translated from Russian), {\it Encycl. Math. Sci.}, vol. 10,
  Several Complex Variables IV, 1990, 195--244, Springer-Verlag.

\bibitem[B]{B} Beloshapka, V. K., On the dimension of the group of
automorphisms of an analytic hypersurface (translated from Russian),
{\it Math. USSR-Izv.} 14(1980), 223--245.

\bibitem[C]{C} Cartan, \'E., Sur les domaines born\'es homog\`enes de l'espace
de $n$ variables complexes, {\it Abh. Math. Sem. Univ. Hamburg} 11(1935),
116--162.

\bibitem[CM]{CM} Chern, S. S. and Moser, J. K., Real hypersurfaces in complex
manifolds, {\it Acta Math.} 133(1974), 219--271.

\bibitem[DY]{DY} Dadok, J. and Yang, P., Automorphisms of tube domains and
spherical hypersurfaces, {\it Amer. J. Math.} 107(1985), 999-1013.

\bibitem[DKR]{DKR} Doubrov, B. M., Komrakov, B. P., and Rabinovich, M.,
Homogeneous surfaces in the three-dimensional affine geometry, in: F.~Dillen
{\it et al.} (eds.), {\it Geometry and Topology of Submanifolds, VIII}, World
Scientific, Singapore, 1996, 168--178.

\bibitem[D]{D} Dynkin, E. B., Maximal subgroups of the classical groups (translated from Russian), {\it Amer. Math. Soc. Transl., II}, Ser. 6(1957), 245--378.

\bibitem[EE1]{EE1} Eastwood, M. G. and Ezhov, V. V., On affine normal forms and
a classification of homogeneous surfaces in affine three-space,
{\it Geom. Dedicata} 77(1999), 11--69.

\bibitem[EE2]{EE2} Eastwood, M. G. and Ezhov, V. V., A classification of
non-degenerate homogeneous equiaffine hypersurfaces in four complex dimensions,
{\it Asian J. Math.} 5(2001), 721--740.

\bibitem[EE3]{EE3} Eastwood, M. G. and Ezhov, V. V., Homogeneous hypersurfaces
with isotropy in affine four-space.
{\it Tr. Mat. Inst. Steklova} 235(2001), 57--70.

\bibitem[E]{E} Ezhov, V. V., Linearization of stability group for a class of
hypersurfaces (translated from Russian), {\it Russian Math. Surveys}
41(1986),203--204.

\bibitem[G]{G} Gindikin, S. G. (ed.), {\it Topics in Geometry: in Memory of
Joseph D'Atri}, Prog. Nonlinear Diff. Eqs. Appls. vol. 20, Birkh\"auser,
Boston, 1996.

\bibitem[I1]{I1} Isaev, A. V., Classification of spherical tube hypersurfaces
that have two minuses in the Levi signature form
(translated from Russian), {\it Math. Notes} 46(1989), 517--523.

\bibitem[I2]{I2} Isaev, A. V., Global properties of spherical tube
hypersurfaces, {\it Indiana Univ. Math. J.} 42(1993), 179--213.

\bibitem[IK]{IK} Isaev, A. V. and Krantz, S. G., On the automorphism groups of hyperbolic manifolds, {\it J. Reine Angew. Math.} 534(2001), 187--194. 

\bibitem[IM]{IM} Isaev, A. V. and Mishchenko, M. A., Classification of
spherical tube hypersurfaces that have one minus in the Levi signature form
(translated from Russian),  {\it Math. USSR-Izv.} 33(1989), 441--472.

\bibitem[K]{K} Komrakov, B. P., Maximal subalgebras of real Lie algebras and a
problem of Sophus Lie (translated from Russian), {\it Soviet Math. Dokl.}
41(1990), 269--273.

\bibitem[L]{L} Loboda, A. V., On local automorphisms of real-analytic
hypersurfaces (translated from Russian), {\it Math. USSR-Izv.} 20(1983),
27--33.

\bibitem[N]{N} Nakajima, K., Homogeneous hyperbolic manifolds and homogeneous Siegel domains, {\it J. Math. Kyoto Univ.} 25(1985), 269--291.

\bibitem[NS1]{NS1} Nomizu, K. and Sasaki, T., A new model of
unimodular-affinely homogeneous surfaces, {\it Manuscripta Math.} 73(1991),
39--44.

\bibitem[NS2]{NS2} Nomizu, K. and Sasaki, T., {\it Affine Differential
Geometry}, Cambridge University Press, Cambridge, 1994.

\bibitem[P-S]{P-S} Pyatetskii-Shapiro, I., {\it Automorphic Functions and the
Geometry of Classical Domains (translated from Russian)}, Gordon and Breach,
1969.

\bibitem[R]{R} Rudin, W., {\it Function Theory in the Unit Ball of $\CC^n$}, Springer-Vergal, 1980.

\bibitem[T]{T} Tanaka, N., On generalized graded Lie algebras and geometric
structures, {\it J. Math. Soc. Japan} 19(1967), 215--254.

\end{thebibliography}
\end{document}